\documentclass{amsart}
\newtheorem{theorem}{Theorem}[section]
\newtheorem{lemma}[theorem]{Lemma}

\theoremstyle{definition}

\newtheorem{ex}[theorem]{Example}

\newtheorem{conj}[theorem]{Conjecture}
\newtheorem{pro}{Proposition}[section]

\theoremstyle{remark}

\numberwithin{equation}{section}

\usepackage[colorlinks]{hyperref}

\begin{document}
	\title{New Results on Generalization of Jordan Centralizers over Matrix Rings}

	\author{Arindam Ghosh}
	\address{Department of Mathematics, Government Polytechnic Kishanganj, Kishanganj-855116}
	\curraddr{}
	\email{E-mail: arindam.rkmrc@gmail.com}
	\thanks{}
	
	\author{Om Prakash$^{\star}$}
	\address{Department of Mathematics, Indian Institute of Technology Patna, Patna-801106}
	\curraddr{}
	\email{om@iitp.ac.in}
	\thanks{* Corresponding author}
	
    	\author{Sushma Singh}
	\address{Department of Mathematics, School of Applied Sciences, KIIT Bhubaneswar- 751 024}
	\curraddr{}
	\email{E-mail: sushma.singhfma@kiit.ac.in}
	\thanks{}
	
	
	
	\date{}
	
	\dedicatory{}
	
	\maketitle
\begin{abstract}
This paper presents a study on Jordan maps over matrix rings with some functional equations related to additive maps on these rings. We first show that every Jordan left (right) centralizer over a matrix ring is a left (right) centralizer. Moreover, every two-sided centralizer over the matrix ring is of a particular form. Further, we prove that any additive map satisfying functional equations over matrix rings becomes a two-sided centralizer. Finally, we conclude our work with some results on the Jordan left $\star$- centralizer over matrix rings and establish some results on functional equations that arise for the $\star$-centralizer.	
\end{abstract}

\section{Introduction}

Throughout, $R$ represents an associative ring and $Z(R)$ is its centre. Recall that a ring $R$ is said to be \emph{prime} if $xRy=0$ for some $x,y\in R$ implies either $x=0$ or $y=0$ and \emph{semiprime} if $xRx=0$ for some $x\in R$ implies $x=0$. The ring $R$ is called \emph{$n$-torsion free} if there exists $x\in R$ such that $nx=0$ implies $x=0$, where $n\geq 2$ is an integer. A map $T:R\rightarrow R$ is said to be a \emph{left (right) centralizer} if $T(x+y)=T(x)+T(y)$ and $T(x y)=T(x)y$ ($T(x y)=xT(y)$), for all $x,y\in R$. It is well known that if $R$ has an identity element $1\neq 0$ and $T:R\rightarrow R$ is a left (right) centralizer, then $T(x)=T(1)x~(T(x)=xT(1))$, for all $x\in R$. The map $T$ is \emph{two-sided centralizer} if it is additive and $T(x y)=T(x)y=xT(y)$, for all $x,y\in R$. Also, the map $T:R\rightarrow R$ is said to be a \emph{Jordan left (right) centralizer} if it is additive and $T(x^2)=T(x)x$ ($T(x^2)=xT(x)$), for all $x\in R$.\\

We can easily see that every left centralizer over a ring is a Jordan left centralizer, but the converse need not be valid (Example \ref{ex2.1}). In 1952, Wendel \cite{wendel1952left}, introduced left centralizer over complex group algebra. After that, Johnson \cite{johnson1964introduction} formally introduced left (right) centralizer over an associative semi-group in 1964. Jordan left centralizer over some rings to be the left centralizer was seen in the past few years by researchers. Interestingly, every left centralizer over a ring is a Jordan left centralizer, but the converse need not be true. Later, in 1992, Bre\v{s}ar and Zalar \cite{brevsar1992structure} proved any Jordan left (right) centralizer over prime rings of characteristic not equal to two is a left (right) centralizer. Also, Zalar proved the same result for semiprime rings \cite{zalar1991centralizers}. Motivated by the above result of Zalar, we show every Jordan left (right) centralizer over any matrix ring is a left (right) centralizer. The importance of the work is that many matrix rings are not semiprime rings. Further, {\cite[Theorem 2.3.2]{beidar1995rings}} motivates to prove that every two-sided centralizer over matrix ring is of a particular form.
In 2008 and 2010, Vukman \cite{vukman2008m,vukman2010m}, introduced $(m, n)$-Jordan derivation and $(m,n)$-Jordan centralizer, respectively. Few more works on $(m, n)$-Jordan derivation and $(m,n)$-Jordan centralizer are available in \cite{fovsner2013note,ali2014generalized,ghosh2019new}.
On the other hand, in 1999, Vukman \cite{vukman1999identity} proved that any additive map $T$ over $2$-torsion-free semiprime ring $R$ with the condition $2T(x^2) = T(x)x + xT(x)$ for all $x \in R$, becomes a two-sided centralizer. Under the same condition, we prove the result for the matrix ring over any $ 2$-torsion free ring. In 2001, Vukman \cite{vukman2001centralizers} proved that an additive map $T$ on a $2$-torsion free semiprime ring $R$ with $T(xyx)=xT(y)x$, for all $x,y\in R$, becomes a two-sided centralizer. We prove the same result for matrix ring $M_r(R)$ ($r\geq 2$ is an integer) over an arbitrary ring $R$.
Again, in 2003, Vukman and Ulbl \cite{vukman2003centralizers} proved that an additive map $T$ on a $2$-torsion free semiprime ring $R$ with $2T(xyx) = T(x)yx + xyT(x)$, for all $x,y\in R$, becomes a two-sided centralizer. We prove the result for matrix ring $M_r(R)$ over the ring $R$.
In the same year, Vukman and Ulbl \cite{vukman2003equation} proved that an additive map $T$ on a $2$-torsion free semiprime ring $R$ with $3T(xyx) = T(x)yx+xT(y)x+xyT(x)$, for all $x,y\in R$, becomes a two-sided centralizer. We prove the result for matrix ring $M_r(R)$ over $2$-torsion free ring $R$. Note that the result is not true for $2$-torsion rings.\\

An involution $\star$ over a ring $R$ is an additive map satisfying $(xy)^{\star}=y^{\star}x^{\star}$ and $(x^{\star})^{\star}=x$, for all $x,y \in R$.
Also, an additive map $T:R\rightarrow R$ is a left (right) $\star$- centralizer if $T(xy)=T(x)y^{\star}$ ($T(xy)=x^{\star}T(y)$), for all $x,y\in R$.
An additive map $T:R\rightarrow R$ is said to be a $\star$- centralizer if $T(xy)=T(x)y^{\star} =x^{\star}T(y)$, for all $x,y\in R$.
An additive map $T:R\rightarrow R$ is said to be a Jordan left (right) $\star$- centralizer if $T(x^2)=T(x)x^{\star}$ ($T(x^2)=x^{\star}T(x)$), for all $x\in R$.
An additive map $T:R\rightarrow R$ is said to be a reverse left (right) $\star$- centralizer if $T(xy)=T(y)x^{\star}$ ($T(xy)=y^{\star}T(x)$), for all $x,y\in R$.
An additive map $T:R\rightarrow R$ is a reverse $\star$-centralizer if $T(xy)=T(y)x^{\star}=y^{\star}T(x)$, for all $x,y\in R$.
In 2013, Ali et al. \cite{ali2013jordan} proved that every Jordan left $\star$-centralizer on a semiprime ring with involution $\star$ and of characteristic different from two
is a reverse left $\star$-centralizer. We prove some results based on the Jordan left $\star$-centralizer over matrix rings and some function equations arising from $\star$-centralizer.

\section{Jordan centralizers over matrix rings}
Let $R$ be a ring with unity $1\neq 0$, $M_r(R)$, $r\geq 2$ be the ring of $r\times r$ matrices over $R$ and $e_{ij}$ be the $r\times r$ matrix with $1$ at $(i,j)$-th place and $0$ elsewhere.
We know that every left centralizer is a Jordan left centralizer, but the converse is invalid. Towards this, we have the following example.
\begin{ex}
\label{ex2.1}
Let $\mathbb{R}$ be the field of real numbers and $S=M_n(\mathbb{R})$. Then $F=S \times S \times S$ is a ring under componentwise addition and for any $(x_1,y_1,z_1),(x_2,y_2,z_2)\in F$, multiplication is defined by $$(x_1,y_1,z_1).(x_2,y_2,z_2)=(0,0,x_{1}y_{2}-x_2y_1).$$ Note that $X^2=0$ for all $X\in F$. Let $P=(1,0,0)$ and $Q=(0,1,0)$. Then $PQ=(0,0,1)$. \\
Suppose $R^{\prime}=\left \{\left( \begin{array}{ccc}
0 & A & B \\
0 & 0 & A\\
0 & 0 & 0 \end{array} \right)\vert \enspace A, B \in F \right \} $ and define $T:R^{\prime} \rightarrow R^{\prime}$ by $$T \left( \begin{array}{ccc}
0 & A & B \\
0 & 0 & A\\
0 & 0 & 0 \end{array} \right)  = \left( \begin{array}{ccc}
0 & 0 & B \\
0 & 0 & 0\\
0 & 0 & 0 \end{array} \right).$$ It can be easily proved that $T$ is a Jordan left centralizer.
Now, consider
\begin{center}
$\tilde{A}=\left( \begin{array}{ccc}
0 & P & 0 \\
0 & 0 & P\\
0 & 0 & 0 \end{array} \right)$ and $\tilde{B}=\left( \begin{array}{ccc}
0 & Q & 0 \\
0 & 0 & Q\\
0 & 0 & 0 \end{array} \right).$
\end{center}
Then $T(\tilde{A}\tilde{B}) \neq T(\tilde{A})\tilde{B}$. Hence, $T$ is not a left centralizer.
\end{ex}

\begin{pro}
If $R$ is a ring and $T:R\rightarrow R$ is a Jordan left centralizer, then $T(xy+yx)=T(x)y+T(y)x$ for all $x,y\in R$.
\end{pro}

\begin{proof}
Substituting $x+y$ for $x$ in $T(x^2)=T(x)x$, we get the result.
\end{proof}

We frequently use this proposition in the proof of Theorem \ref{lc1}.
\begin{theorem}
\label{lc1}
Let $R$ be a ring. Then every Jordan left (right) centralizer $T: M_r(R)\rightarrow M_r(R)$ is a left (right) centralizer.
\end{theorem}
\begin{proof}
Let $T$ be a Jordan left centralizer on $M_r(R)$ and for all $i,j\in \{1,2,\dots ,r\}$,
\begin{equation}
\label{eq:mnr1}
T(e_{ij})=\sum_{k=1}^{r}\sum_{l=1}^{r}a_{kl}^{(ij)}e_{kl}~, ~\text{for} ~a_{kl}^{(ij)} \in R.
\end{equation}

Since $e_{ii}^2=e_{ii}$ and $T$ is a Jordan left centralizer on $M_r(R)$, we have
\begin{equation}
\label{eq:mnr2}
T(e_{ii})=\sum_{k=1}^{r}a_{ki}^{(ii)}e_{ki}~, ~\text{for all} ~i \in \{1,2,\dots ,r\}.
\end{equation}

Also, $e_{ij}=e_{ii}e_{ij}+e_{ij}e_{ii}$, for $i\neq j$. Therefore,

\begin{equation}
\label{eq:mnr3}
T(e_{ij})=\sum_{k=1}^{r}a_{ki}^{(ii)}e_{kj}+{\sum_{k=1}^{r}a_{ki}^{(ij)}e_{kj}}.
\end{equation}

Again, we have $e_{ij}=e_{ij}e_{jj}+e_{jj}e_{ij}$. Hence, by applying \eqref{eq:mnr2} and \eqref{eq:mnr3}, we have
\begin{equation}
\label{eq:mnr4}
T(e_{ij})=\sum_{k=1}^{r}a_{ki}^{(ii)}e_{kj}.
\end{equation}

Now, let $s\in R$ and for all $i,j\in \{1,2,\dots ,r\}$,
\begin{equation}
\label{eq:mnr5}
T(se_{ij})=\sum_{k=1}^{r}\sum_{l=1}^{r}a_{kl}^{s(ij)}e_{kl}~, ~\text{for} ~a_{kl}^{s(ij)} \in R.
\end{equation}

We know $se_{ij}=(se_{ij})e_{jj}+e_{jj}(se_{ij})$, for $i\neq j$. Since, $T$ is a Jordan left centralizer on $M_r(R)$, applying \eqref{eq:mnr2} and \eqref{eq:mnr5},
\begin{equation}
\label{eq:mnr6}
T(se_{ij})=\sum_{k=1}^{r}a_{kj}^{s(ij)}e_{kj}.
\end{equation}

Similarly, by $se_{ij}=e_{ii}(se_{ij})+(se_{ij})e_{ii}$, \eqref{eq:mnr2} and \eqref{eq:mnr6}, we have
\begin{equation}
\label{eq:mnr7}
T(se_{ij})=\sum_{k=1}^{r}a_{ki}^{(ii)}se_{kj}.
\end{equation}

Now, $2se_{ii}=e_{ii}(se_{ii})+(se_{ii})e_{ii}$, so by \eqref{eq:mnr2} and \eqref{eq:mnr5}, we have
\begin{equation}
\label{eq:mnr8.1}
\begin{aligned}
& 2T(se_{ii})= T(e_{ii})(se_{ii})+T(se_{ii})e_{ii}\implies a_{ki}^{s(ii)}=a_{ki}^{(ii)}s, ~\text{for all}~k=1,2,
\dots,r.
\end{aligned}
\end{equation}

For all $i \neq j$, $(se_{ii})(e_{jj})+(e_{jj})(se_{ii})=0$, by \eqref{eq:mnr2}, we have
\begin{equation}
\label{eq:mnr8.2}
\begin{aligned}
& 0=T(se_{ii})e_{jj}+T(e_{jj})se_{ii}\implies a_{kj}^{s(ii)}=0, ~\text{for all}~k=1,2,
\dots,r.
\end{aligned}
\end{equation}

From \eqref{eq:mnr8.1} and \eqref{eq:mnr8.2},
\begin{equation}
\label{eq:mnr8}
\begin{aligned}
 T(se_{ii})=\sum_{k=1}^{r}a_{ki}^{(ii)}se_{ki}.
\end{aligned}
\end{equation}

By the relation \eqref{eq:mnr2}, \eqref{eq:mnr4}, \eqref{eq:mnr7} and \eqref{eq:mnr8}, we have $T(se_{ij})=T(e_{ij})s$, for all $s\in R$.

Now, by \eqref{eq:mnr2}, $T(1)=\sum_{i=1}^{r}\sum_{k=1}^{r}a_{ki}^{(ii)}e_{ki}$.

Let $A=T(1)$ and $X=\sum_{i=1}^{r}\sum_{j=1}^{r}x_{ij}e_{ij}$ be any element in $M_r(R)$ (Here $1$ denotes the identity matrix in $M_r(R)$).
\begin{align*}
T(X)&=T \left(\sum_{i=1}^{r}\sum_{j=1}^{r}x_{ij}e_{ij}\right)
=\sum_{i=1}^{r}\sum_{j=1}^{r}T(e_{ij})x_{ij}\\
&=\sum_{i=1}^{r}\sum_{j=1}^{r}\left(\sum_{k=1}^{r}a_{ki}^{(ii)}e_{kj}\right)x_{ij}\\
&=\left(\sum_{i=1}^{r}\sum_{k=1}^{r}a_{ki}^{(ii)}e_{ki}\right)\left(\sum_{i=1}^{r}\sum_{j=1}^{r}x_{ij}e_{ij}\right)
=AX.
\end{align*}
Hence $T$ is a left centralizer.

Due to symmetry, every Jordan right centralizer over $M_r(R)$ is a right centralizer.
\end{proof}

\begin{lemma}
\label{mnr1}
Let $R$ be a ring and $T: M_r(R)\rightarrow M_r(R)$ be a two-sided centralizer. Then there exists an $\alpha\in Z(R)$ such that $T(X)=\alpha X$, for all $X\in M_r(R)$.
\end{lemma}

\begin{proof} Let $1$ be the identity matrix in $M_r(R)$.
Since $T$ is a two-sided centralizer, $T(X)=T(1)X=XT(1)$, for all $X\in M_r(R)$. Then $T(1)$ commutes with every element of $M_r(R)$. It is well-known that the center of a matrix ring coincides with the center of a ring (embedded diagonally into the matrix ring). Hence, $T(X)=\alpha X$ for all $X\in M_r(R)$ where $T(1)=\alpha .1$.
\end{proof}

\begin{theorem}
\label{lc2}
Let $m \geq 1$, $n \geq 1$ and  $m,~n\in \mathbb{Z}$, $R$ be a ring with $n(m+n)^3$-torsion free. If $T: M_r(R) \rightarrow M_r(R)$ be an additive mapping such that there exists a two-sided centralizer $T_0: M_r(R) \rightarrow M_r(R)$ satisfying
\begin{equation}
\label{eq:mnr12}
(m  + n)T(x^2) = mT(x)x + nxT_0(x),~\text{for all}~ x \in M_r(R),
\end{equation}
then $T$ becomes a two-sided centralizer. In fact, $T=T_0$.
\end{theorem}

\begin{proof}
Since $T$ satisfies \eqref{eq:mnr12} and $T_0$ is a two sided centralizer on $M_r(R)$, by Lemma \ref{mnr1}, we have
\begin{equation}
\label{eq:mnr13}
(m  + n)T(x^2) = mT(x)x + n\alpha x^2,~\text{for all}~ x \in M_r(R)~\text{for some}~ \alpha \in Z(R).
\end{equation}

Replacing $x$ by $x+y$ in \eqref{eq:mnr13},
\begin{equation}
\label{eq:mnr14}
(m  + n)T(xy+yx) =  n\alpha (xy+yx)+mT(x)y+mT(y)x,~\text{for all}~ x,y \in M_r(R).
\end{equation}

Let $T(e_{ij})$ be of the form \eqref{eq:mnr1}. Since $e_{ii}^2=e_{ii},$ using \eqref{eq:mnr1} and \eqref{eq:mnr13},
\begin{equation}
\label{eq:mnr15}
(m  + n)T(e_{ii}) = m\underset{k\neq i}{\sum_{k=1}^{r}}a_{ki}^{(ii)}e_{ki}+(ma_{ii}^{(ii)}+n\alpha)e_{ii}.
\end{equation}

Now, $e_{ij}=e_{ii}e_{ij}+e_{ij}e_{ii}$, for $i\neq j$. By \eqref{eq:mnr14} and  \eqref{eq:mnr15},
\begin{equation}
\label{eq:mnr16}
\begin{aligned}
&(m  + n)^2T(e_{ij})\\
& = m^2\underset{k\neq i}{\sum_{k=1}^{r}}a_{ki}^{(ii)}e_{kj}+m(ma_{ii}^{(ii)}+n\alpha)e_{ij} + m(m+n){\sum_{k=1}^{r}}a_{ki}^{(ij)}e_{ki}+(m+n)n\alpha e_{ij}.
\end{aligned}
\end{equation}

Also, we have $e_{ij}=e_{ij}e_{jj}+e_{jj}e_{ij}$. Applying \eqref{eq:mnr15} and \eqref{eq:mnr16},
\begin{equation}
\label{eq:mnr17}
\begin{aligned}
&(m  + n)^3T(e_{ij})\\
& = m^3\underset{k\neq i}{\sum_{k=1}^{r}}a_{ki}^{(ii)}e_{kj}+m^2(ma_{ii}^{(ii)}+n\alpha)e_{ij} + m(m+n)n\alpha e_{ij}+(m+n)^2n\alpha e_{ij}.
\end{aligned}
\end{equation}

Let $1\neq 0$ be the identity element in $M_r(R)$. Since $1^2=1$, $1=\sum_{k=1}^{r}e_{kk}$, and $R$ is $n$-torsion free,
\begin{equation}
\label{eq:mnr18}
\begin{aligned}
(m  + n)\sum_{k=1}^{r}T(e_{kk})=(m  + n)\alpha \left(\sum_{k=1}^{r}e_{kk}\right).
\end{aligned}
\end{equation}

Applying \eqref{eq:mnr15} to \eqref{eq:mnr18}, we have
\begin{equation}
\label{eq:mnr19}
\begin{aligned}
a_{ij}^{(kk)}=0~\text{and}~a_{kk}^{(kk)}=\alpha, \forall~i,j,k~\text{with}~i\neq j~(\text {Since $R$ is $(m+n)$-trosion free}).
\end{aligned}
\end{equation}

Again, by \eqref{eq:mnr19}, \eqref{eq:mnr15} and \eqref{eq:mnr17},
\begin{equation}
\label{eq:mnr20}
\begin{aligned}
T(e_{ij})=\alpha e_{ij},~\text{for any}~i,j~(\text {Since $R$ is $(m+n)^3$-trosion free}).
\end{aligned}
\end{equation}

Therefore, by \eqref{eq:mnr20}, $T=T_0$,  and hence $T$ is a two-sided centralizer.
\end{proof}

\section{Some functional equations over matrix rings}

The primary purpose of this section is to investigate a map satisfying some equations to become a two-sided centralizer.

\begin{theorem}
\label{thm2.3}
Let $R$ be a ring with $2$-torsion free. If $T$ is an additive function on $M_r(R)$ which satisfies
\begin{equation}
\label{eq:eq1}
2T(x^2)=T(x)x+xT(x), ~\text{for all} ~x\in M_r(R),
\end{equation}
then $T$ become a two-sided centralizer. In particular, $T( x)=\alpha x$, for all $x\in M_r(R)$ and for an $\alpha \in Z(R)$.
\end{theorem}

\begin{proof}
Let $T(e_{ij})$ and $T(se_{ij})$ be of the form \eqref{eq:mnr1} and \eqref{eq:mnr5}, respectively for $s\in R$. Since $e_{ii}^2=e_{ii}$, $R$ is $2$-torsion free and satisfies \eqref{eq:eq1}, we have
\begin{equation}
\label{eq:eq2}
T(e_{ii})=a_{ii}^{(ii)}e_{ii}~, ~\text{for all} ~i \in \{1,2,\dots ,r\}.
\end{equation}

Substituting $x+y$ for $x$ in \eqref{eq:eq1},
\begin{equation}
\label{eq:eq3}
2T(xy+yx)=T(x)y+T(y)x+xT(y)+yT(x)~, ~\text{for all} ~x,y \in M_r(R).
\end{equation}

For $i\neq j$, applying \eqref{eq:eq3} on $e_{ij}=e_{ii}e_{ij}+e_{ij}e_{ii}$, all the other coefficients of $T(e_{ij})$ except $e_{ij}$ become zero and
\begin{equation}
\label{eq:eq4}
a_{ij}^{(ij)}=a_{ii}^{(ii)}, ~\text{for all} ~i,j \in \{1,2,\dots ,r\}.
\end{equation}

Hence, we have
\begin{equation}
\label{eq:eq5}
T(e_{ij})=a_{ij}^{(ij)}e_{ij}=a_{ii}^{(ii)}e_{ij}, ~\text{for all} ~i,j \in \{1,2,\dots ,r\}.
\end{equation}

Again, applying \eqref{eq:eq3} on $e_{ij}=e_{ij}e_{jj}+e_{jj}e_{ij}$, we have
\begin{equation}
\label{eq:eq6}
a_{ii}^{(ii)}=a_{jj}^{(jj)}, ~\text{for all} ~i,j \in \{1,2,\dots ,r\}.
\end{equation}

Also,
\begin{equation}
\label{eq:eq7}
T(e_{ij})=a_{11}^{(11)}e_{ij}, ~\forall ~i,j \in \{1,2,\dots ,r\}.
\end{equation}

For $i\neq j$ and $s\in R$, applying \eqref{eq:eq3} on $se_{ij}=(se_{ij})e_{jj}+e_{jj}(se_{ij})$, all the other coefficients of $T(se_{ij})$ except $e_{ij}$ become zero and
\begin{equation}
\label{eq:eq8}
a_{ij}^{s(ij)}=sa_{11}^{(11)}, ~\forall ~i,~j .
\end{equation}

For $i\neq j$ and $s\in R$,
\begin{equation}
\label{eq:eq9}
T(se_{ij})=sa_{11}^{(11)}e_{ij}, ~\text{for all} ~i,~ j .
\end{equation}

Applying \eqref{eq:eq3}, $2se_{ii}=(se_{ii})e_{ii}+e_{ii}(se_{ii})$,
\begin{equation}
\label{eq:eq10}
2a_{ii}^{s(ii)} = a_{11}^{(11)}s+sa_{11}^{(11)}, ~\text{for all} ~i \in \{1,2,\dots ,r\}.
\end{equation}

For all $i \neq j$, $(se_{ii})(e_{jj})+(e_{jj})(se_{ii})=0$, by \eqref{eq:eq3},
\begin{equation}
\label{eq:eq10.1}
\begin{aligned}
a_{kj}^{s(ii)}=0~\text{and}~a_{jk}^{s(ii)}=0, ~\text{for all}~k=1,2,
\dots,r~(\text{Since} ~R~ \text{is}~\text{2-torsion free}).
\end{aligned}
\end{equation}

Applying \eqref{eq:eq3}, $se_{ij}=(se_{ij})e_{ii}+e_{ii}(se_{ij})$,
\begin{equation}
\label{eq:eq11}
a_{11}^{(11)}s=sa_{11}^{(11)}.
\end{equation}

Again, we have
\begin{equation}
\label{eq:eq12}
T(se_{ii})=sa_{11}^{(11)}e_{ii}, ~\text{for all} ~i \in \{1,2,\dots ,r\}.
\end{equation}

Let $X=\sum_{i=1}^{r}\sum_{j=1}^{r}x_{ij}e_{ij}$, $x_{ij} \in R$. Let $\alpha=a_{11}^{(11)}$. Then by using \eqref{eq:eq9} and \eqref{eq:eq12}, we get
\begin{equation}
\label{eq:eq13}
T(x)=a_{11}^{(11)}x=\alpha x, ~\text{for all}~x\in M_r(R).
\end{equation}

Thus, by \eqref{eq:eq11} and \eqref{eq:eq13}, $T$ is a two-sided centralizer.
\end{proof}

\begin{ex}
\label{rem1}
Let $\mathbb{Z}_2$ be the ring of integers modulo $2$.
Now, we provide an example of an additive map $T$ on $M_r(R)$ with $2T(X^2)=T(X)X+XT(X)$, for all $X \in M_r(R)$, but it is not a two-sided centralizer.

Let $X= \begin{bmatrix}
x & y \\
z & t \end{bmatrix}\in M_2(\mathbb{Z}_2)$ and $T:M_2(\mathbb{Z}_2)\rightarrow M_2(\mathbb{Z}_2)$ is defined by
\begin{center}
$T(X)= \begin{bmatrix}
x+y+z+t & 0 \\
0 & x+y+z+t \end{bmatrix}$.
 \end{center}

 Then $T$ satisfies $2T(X^2)=T(X)X+XT(X)$, for all $X \in M_2(\mathbb{Z}_2)$. Also, for $X=e_{11}$ and $Y= e_{12}$, $T(XY)\neq T(X)Y$. Thus, $T$ is not a two-sided centralizer.
\end{ex}

\begin{theorem}
\label{thm2.4}
Let $R$ be a $ 2$-torsion free ring. If $T$ is an additive mapping on $M_r(R)$ satisfying
\begin{equation}
\label{eq:eq14}
2T(x^2)=T(x)x+xT_0(x), ~\text{for all} ~x\in M_r(R),
\end{equation}
where $T_0$ is an additive mapping on $M_r(R)$ satisfying \begin{equation}
\label{eq:eq15}
2T_0(x^2)=T_0(x)x+xT_0(x), ~\text{for all} ~x\in M_r(R),
\end{equation}
then $T$ is a two-sided centralizer.
\end{theorem}

\begin{proof}
By Theorem \ref{thm2.3} and Lemma \ref{mnr1}, $T_0$ is a two-sided centralizer and $T_0(x)=\alpha x$ for some $\alpha \in Z(R)$ and $x\in M_r(R)$. Therefore, we get
\begin{equation}
\label{eq:eq16}
2T(x^2)=T(x)x+\alpha x^2, ~\text{for all} ~x\in M_r(R).
\end{equation}

Again, let $T(e_{ij})$ and $T(se_{ij})$ be of the form \eqref{eq:mnr1} and \eqref{eq:mnr5}, respectively for $s\in R$. Applying \eqref{eq:eq16} on $e_{ii}^2=e_{ii}$, all the coefficients of $T(e_{ii})$ except $e_{ii}$ become zero and
\begin{equation}
\label{eq:eq17}
a_{ii}^{(ii)}=\alpha, ~\text{for all} ~i \in \{1,2,\dots ,r\}.
\end{equation}

We have
\begin{equation}
\label{eq:eq18}
T(e_{ii})=\alpha e_{ii}, ~\text{for all} ~i \in \{1,2,\dots ,r\}.
\end{equation}

Linearizing \eqref{eq:eq16},
\begin{equation}
\label{eq:eq19}
2T(xy+yx)=T(x)y+T(y)x+\alpha (xy+yx), ~\text{for all} ~x, y\in M_r(R).
\end{equation}

For $i\neq j$, applying \eqref{eq:eq19} on $e_{ij}=e_{ii}e_{ij}+e_{ij}e_{ii}$, each coefficient of $T(e_{ij})$ other than $e_{ij}$ is zero and
\begin{equation}
\label{eq:eq20}
a_{ij}^{(ij)}=\alpha, ~\text{for all} ~i,j \in \{1,2,\dots ,r\}.
\end{equation}

For $i\neq j$,
\begin{equation}
\label{eq:eq21}
T(e_{ij})=\alpha e_{ij}, ~\forall ~i,j .
\end{equation}

For $i\neq j$ and $s\in R$, applying \eqref{eq:eq19} on $se_{ij}=(se_{ij})e_{jj}+e_{jj}(se_{ij})$, all the coefficients of $T(se_{ij})$ except $e_{ij}$ are zero and
\begin{equation}
\label{eq:eq22}
a_{ij}^{s(ij)}=\alpha s, ~\text{for all} ~i,j \in \{1,2,\dots ,r\}.
\end{equation}

For $i\neq j$ and $s\in R$,
\begin{equation}
\label{eq:eq23}
T(se_{ij})=\alpha s e_{ij}, ~\text{for all} ~i,j \in \{1,2,\dots ,r\}.
\end{equation}

Applying \eqref{eq:eq19} on $2se_{ii}=(se_{ii})e_{ii}+e_{ii}(se_{ii})$, each coefficient of $T(se_{ii}),$ except $e_{ii},$ is zero and
\begin{equation}
\label{eq:eq24}
a_{ii}^{s(ii)}=\alpha s, ~\text{for all} ~i \in \{1,2,\dots ,r\}.
\end{equation}

For all $s\in R$,
\begin{equation}
\label{eq:eq25}
T(se_{ii})=\alpha s e_{ii}, ~\text{for all} ~i \in \{1,2,\dots ,r\}.
\end{equation}

By \eqref{eq:eq23} and \eqref{eq:eq25}, we conclude that
\begin{equation}
\label{eq:eq26}
T=T_0.
\end{equation}

Thus, $T$ is a two-sided centralizer.
\end{proof}

Now, we give an example which shows that $M_r(R)$ is not always a semiprime ring. Hence, Theorem \ref{thm2.3} is not a consequence of any result of Vukman \cite{vukman1999identity}.
\begin{ex}
\label{ex2.9}
Let $\mathbb{Z}_9$ be the ring of residue classes of integers modulo $9$. Then $\mathbb{Z}_9$ is $2$-torsion free with unity $1\neq 0$. It is easy to compute that $(3e_{11})M_2(\mathbb{Z}_9)(3e_{11})=0$. But $3e_{11}\neq 0$. Therefore, $M_2(\mathbb{Z}_9)$ is not a semiprime ring (since in a semiprime ring $R$, $aRa=0 \implies a=0$).
\end{ex}

\begin{theorem}
\label{thm2.10}
Let $R$ be a ring. If $T$ is an additive map on $M_r(R)$ with
\begin{equation}
\label{eq:eq27}
T(xyx)=xT(y)x, ~\text{for all} ~x,y\in M_r(R),
\end{equation}
then $T$ becomes a two-sided centralizer. In particular, $T( x)=\alpha x$ for all $x\in M_r(R)$ and for an $\alpha \in Z(R)$.
\end{theorem}

\begin{proof}
Put $x=0,~y=1$ ($1$ and $0$ denote the identity and zero matrix, respectively) in \eqref{eq:eq27}, we have $T(0)=0$. Let
$T(e_{ij})$ and $T(se_{ij})$ be of the form \eqref{eq:mnr1} and \eqref{eq:mnr5}, respectively for $s\in R$. Since $e_{ii}^3=e_{ii}$ and $x.1.x=x^2$, using \eqref{eq:eq27},
\begin{equation}
\label{eq27a}
\begin{aligned}
&T(e_{ii})=a_{ii}^{(ii)}e_{ii},\\
&(\text{hence})~T(1)=\sum_{k=1}^{r}a_{kk}^{(kk)}e_{kk},\\
& T(x^2)=xT(1)x ~\text{for all} ~x,y \in M_r(R).
\end{aligned}
\end{equation}

Put $x=x+y$ in \eqref{eq27a},
\begin{equation}
\begin{aligned}
\label{eq27b}
T(xy+yx)=xT(1)y+yT(1)x,~\text{for all} ~x,y \in M_r(R).
\end{aligned}
\end{equation}

Put $x=e_{ii}$ and $y=e_{ij}$ ($i\neq j$) in \eqref{eq27b},
\begin{equation}
\begin{aligned}
\label{eq27c}
& T(e_{ij})=a_{ii}^{(ii)}e_{ij}\\
& (\text{also putting}~x=e_{ij},~y=e_{jj}~\text{in}~\eqref{eq27b}),~ T(e_{ij})=a_{jj}^{(jj)}e_{ij}\\
&(\text{hence})~a_{ii}^{(ii)}=a_{jj}^{(jj)}=\alpha~(\text{say})\\
&T(e_{ij})=\alpha e_{ij}.
\end{aligned}
\end{equation}

Since $se_{ij}=(se_{ij})e_{jj}+e_{jj}(se_{ij})$, using \eqref{eq27c},
\begin{equation}
\begin{aligned}
\label{eq27d}
& T(se_{ij})=s\alpha e_{ij}\\
& (\text{similarly}),~ T(se_{ij})=\alpha s e_{ij}\\
&(\text{hence})~\alpha s=s \alpha~(s\in R).
\end{aligned}
\end{equation}

Put $x=x+z$ in \eqref{eq:eq27},
\begin{equation}
\begin{aligned}
\label{eq27e}
T(xyz+zyx)=xT(y)z+zT(y)x,~\text{for all}~x,y,z\in M_r(R).
\end{aligned}
\end{equation}

Put $x=se_{ii}$, $y=e_{ji}$ and $z=e_{ii}$ in \eqref{eq27e},
\begin{equation}
\label{eq27f}
T(se_{ii})=s\alpha e_{ii}=\alpha se_{ii}.
\end{equation}

From \eqref{eq27d} and \eqref{eq27f}, $T(x)=\alpha x$ for all $x\in M_r(R)$.
\end{proof}

\begin{theorem}
\label{thm2.11}
Let $R$ be a ring. If $T$ is an additive map on $M_r(R)$ with
\begin{equation}
\label{eq:eq28}
2T(xyx) = T(x)yx + xyT(x), ~\text{for all} ~x,y\in M_r(R),
\end{equation}
then $T$ is a two-sided centralizer. In particular, there exists an $\alpha \in Z(R)$ such that $T( x)=\alpha x$, $\forall~x\in M_r(R)$.
\end{theorem}

\begin{proof}
Put $x=0,~y=0$ in \eqref{eq:eq28}, we have $2T(0)=0$. Let
$T(e_{ij})$ and $T(se_{ij})$ be of the form \eqref{eq:mnr1} and \eqref{eq:mnr5}, respectively for $s\in R$. Putting $x=e_{ii}$, $y=e_{ii}$ in \eqref{eq:eq28},
\begin{equation}
a_{il}^{(ii)}=0 ~\text{for all}~l \in \{1,2,\dots,i-1,i+1,\dots,r\}.
\end{equation}

Putting $x=e_{ii}$, $y=e_{ik}$ ($k\neq i$) in \eqref{eq:eq28},
\begin{equation}
a_{kl}^{(ii)}=0 ~\text{for all}~l \in \{1,2,\dots,r\}~\text{and hence}~T(e_{ii})=a_{ii}^{(ii)}e_{ii}.
\end{equation}

Putting $x=e_{ij}$, $y=e_{ji}$ ($j\neq i$) in \eqref{eq:eq28},
\begin{equation}
a_{im}^{(ij)}=0~\text{for all}~m\in \{1,2,\dots,j-1,j+1,\dots,r\}.
\end{equation}

Putting $x=e_{ij}$, $y=e_{jk}$ ($j\neq i$,$i\neq k$) in \eqref{eq:eq28},
\begin{equation}
a_{km}^{(ij)}=0~\text{for all}~m\in \{1,2,\dots,r\} ~\text{and hence}~T(e_{ij})=a_{ij}^{(ij)}e_{ij}.
\end{equation}

Putting $x=1$ and $y=e_{ij}$ in \eqref{eq:eq28}, $2a_{ij}^{(ij)}=a_{ii}^{(ii)}+a_{jj}^{(jj)}$. Also, putting $x=e_{ij}+e_{ii}$ and $y=e_{ji}$ in \eqref{eq:eq28}, $a_{ii}^{(ii)}=a_{ij}^{(ij)}$. Therefore, $a_{ij}^{(ij)}=a_{ii}^{(ii)}=a_{jj}^{(jj)}=\alpha$ (say).
\begin{equation}
\label{eq28a}
 \text{Hence,}~ T(e_{kl})=\alpha e_{kl}~ \text{for all}~ k,~l.
\end{equation}

Putting $x=e_{ii}$ and $y=se_{ii}$ in \eqref{eq:eq28} and using \eqref{eq28a},
\begin{equation}
\begin{aligned}
\label{eq28b}
& 2T(se_{ii})=(\alpha s+s \alpha)e_{ii}\\
& 2a_{ii}^{s(ii)}=\alpha s+s \alpha.
\end{aligned}
\end{equation}

Put $x=x+z$ in \eqref{eq:eq28},
\begin{equation}
\begin{aligned}
\label{eq28c}
2T(xyz+zyx)=T(x)yz+xyT(z)+T(z)yx+zyT(x),~\text{for all}~x,y,z.
\end{aligned}
\end{equation}

Putting $x=se_{ii}$, $y=e_{ij}$ and $z=e_{ji}$ ($j\neq i$) in \eqref{eq28c} and using \eqref{eq28a} and \eqref{eq28b},
\begin{equation}
a_{ii}^{s(ii)}=s \alpha=\alpha s~\text{and}~a_{jl}^{s(ii)}=0 ~\text{for}~l\in \{1,2,\dots,r\}.
\end{equation}

Put $x=e_{ij}$, $y=e_{ji}$ and $z=se_{ii}$ in \eqref{eq28c},
\begin{equation}
a_{il}^{s(ii)}=0~\text{for}~l\in \{1,2,\dots,i-1,i+1,\dots,r\}~\text{and hence}~T(se_{ii})=\alpha s e_{ii}~(s\in R, \alpha\in Z(R)).
\end{equation}

Putting $x=se_{ij}$, $y=e_{jj}$ and $z=e_{jj}$ in \eqref{eq28c},
\begin{equation}
a_{ij}^{s(ij)}=s \alpha ~\text{and}~a_{mj}^{s(ij)}=0~\text{for}~m\in \{1,2,\dots,i-1,i+1,\dots,r\}.
\end{equation}

Putting $x=se_{ij}$, $y=e_{kk}$ and $z=e_{kk}$ ($k\neq j$) in \eqref{eq28c},
\begin{equation}
a_{mk}^{s(ij)}=0~\text{for}~m\in \{1,2,\dots,r\}~\text{and hence}~T(se_{ij})=\alpha se_{ij}.
\end{equation}

Thus, $T(x)=\alpha x$ for all $x\in M_r(R)$.

\end{proof}

\begin{theorem}
\label{thm2.12}
Let $R$ be a $ 2$-torsion free ring. If $T$ is an additive map on $M_r(R)$ with
\begin{equation}
\label{eq:29}
3T(xyx) = T(x)yx+xT(y)x+xyT(x), ~\text{for all} ~x,y\in M_r(R),
\end{equation}
then $T$ is a two-sided centralizer. In particular, there exists an $\alpha \in Z(R)$ such that $T( x)=\alpha x$, for all $x\in M_r(R)$.
\end{theorem}

\begin{proof}
Let $T(e_{ij})$ and $T(se_{ij})$ be of the form \eqref{eq:mnr1} and \eqref{eq:mnr5}, respectively for $s\in R$. Since $e_{ij}e_{ji}e_{ij}=e_{ij}$ and the torsion condition on $R$, using \eqref{eq:29},
\begin{equation}
\label{eq:30}
\begin{aligned}
&a_{1j}^{(ij)}=\dots=a_{i-1,j}^{(ij)}=a_{i+1,j}^{(ij)}=\dots=a_{rj}^{(ij)}=0,\\
&a_{i1}^{(ij)}=\dots=a_{i,j-1}^{(ij)}=a_{i,j+1}^{(ij)}=\dots=a_{ir}^{(ij)}=0,\\
&a_{ij}^{(ij)}=a_{ji}^{(ji)}~\text{for all} ~i,j \in \{1,2,\dots ,r\}.
\end{aligned}
\end{equation}

Substituting $1$ and $0$ for $x$ and $y$ in \eqref{eq:29}, respectively, we have $T(0)=0$. Let $k\neq i$ and $l\neq j$ and since $e_{lk}e_{ij}e_{lk}=0$, using \eqref{eq:29},
\begin{equation}
\label{eq:31}
\begin{aligned}
& a_{kl}^{(ij)}=0.\\
& \text{Hence,}~T(e_{ij})=a_{ij}^{(ij)}e_{ij}.
\end{aligned}
\end{equation}

Linearizing \eqref{eq:29},
\begin{equation}
\label{eq:32}
\begin{aligned}
&3T(xyz+zyx)=T(x)yz+T(z)yx+xT(y)z\\
&+zT(y)x+xyT(z)+zyT(x), ~\text{for all} ~x,y,z \in M_r(R).
\end{aligned}
\end{equation}

Substituting $e_{ii},~e_{ij}~\text{and}~e_{jj}$ for $x,~y$ and $z$ respectively in \eqref{eq:32},
\begin{equation}
\label{eq:33}
\begin{aligned}
2a_{ij}^{(ij)}=a_{ii}^{(ii)}+a_{jj}^{(jj)}~\text{for all}~i\neq j.
\end{aligned}
\end{equation}

Substituting $e_{ji},~e_{ij}~\text{and}~e_{jj}$ for $x,~y$ and $z$,    respectively in \eqref{eq:32},
\begin{equation}
\label{eq:34}
\begin{aligned}
& (\text{by}~\eqref{eq:30})~a_{jj}^{(jj)}=a_{ij}^{(ij)}=a_{ji}^{(ji)}, \\
& (\text{by above and}~\eqref{eq:30})~a_{ii}^{(ii)}= a_{ij}^{(ij)}=a_{ji}^{(ji)}=a_{jj}^{(jj)}, ~\text{for all}~i\neq j.\\
&\text{Hence}~T(e_{ij})=\alpha e_{ij}~\text{(letting,}~ \alpha=a_{11}^{(11)}), ~\text{for all}~i, j.
\end{aligned}
\end{equation}

Let $s\in R$. Substituting $1~\text{and}~se_{ij}$ for $x$ and $y$, respectively in \eqref{eq:29} and using \eqref{eq:34},
\begin{equation}
\label{eq:35}
\begin{aligned}
&2T(se_{ij})=(\alpha s+s \alpha)e_{ij}, \\
& 2a_{ij}^{s(ij)}=\alpha s+s \alpha, \\
&T(se_{ij})=a_{ij}^{s(ij)}e_{ij}, ~\text{for all}~i, j.
\end{aligned}
\end{equation}

Substituting $se_{ii},~e_{ij}~\text{and}~e_{jj}$ for $x,~y$ and $z$, respectively in \eqref{eq:32} and using \eqref{eq:34},
\begin{equation}
\label{eq:36}
\begin{aligned}
a_{ij}^{s(ij)}=s \alpha=\alpha s~(\text{using}~\eqref{eq:35}),~\text{for all}~i\neq j.
\end{aligned}
\end{equation}

Therefore, $T(x)=\alpha x $, for all $x\in M_r(R)$.
\end{proof}

\begin{ex}
\label{ex2.14}
Now, we give an example of an additive mapping $T$ on $M_r(R)$ satisfying \eqref{eq:29}, but it is not a two-sided centralizer. Let $\mathbb{Z}_2$ be the ring of residue classes of integers modulo $2$, which is not $2$-torsion free. Let $X= \begin{bmatrix}
x & y \\
z & t \end{bmatrix}\in M_2(\mathbb{Z}_2)$ and $T:M_2(\mathbb{Z}_2)\rightarrow M_2(\mathbb{Z}_2)$ be defined by
$$ T\left( \begin{bmatrix}
 x & y \\
z & t \end{bmatrix}\right)=
 \begin{bmatrix}
y & 0 \\
x & 0\end{bmatrix}.
$$
Then $T$ satisfies \eqref{eq:29}. In this case, for $X=e_{11}$ and $Y=e_{12}$, $T(XY)=e_{11}\neq e_{22}= T(X)Y$. Thus, $T$ is not a two-sided centralizer.
\end{ex}

\section{Jordan $\star$-centralizers over rings}
It is easy to prove that every reverse left $\star$-centralizer is a Jordan left $\star$-centralizer, but the converse is not generally valid. Hence, it is a genuine attempt to classify rings and algebras over which Jordan left $\star$-centralizer becomes reverse left $\star$-centralizer.

\begin{ex}
Let $S=\mathbb{Z}_2[x,y]$ with $x^2=y^2=0$ and $R_1$ be the subring of $S$ generated by $x$ and $y$. One can easily recognize that every element of $R_1$ is of square zero and $R_1$ is commutative.

Let $R_2=\left \{\left( \begin{array}{ccc}
0 & A & B \\
0 & 0 & A\\
0 & 0 & 0 \end{array} \right)\vert \enspace A, B \in R_1 \right \} $ and  $\left( \begin{array}{ccc}
0 & A & B \\
0 & 0 & A\\
0 & 0 & 0 \end{array} \right)^{\star}=\left( \begin{array}{ccc}
0 & A & B \\
0 & 0 & A\\
0 & 0 & 0 \end{array} \right)$. Then $R_2$ is a ring with involution $\star$.

Define $T:R_2 \rightarrow R_2$ as,~~~~~\\
$T \left( \begin{array}{ccc}
0 & A & B \\
0 & 0 & A\\
0 & 0 & 0 \end{array} \right)=\left( \begin{array}{ccc}
0 & 0 & B \\
0 & 0 & 0\\
0 & 0 & 0 \end{array} \right)$. It is easy to check that $T$ is a Jordan left $\star$-centralizer.

Let $\tilde{X}=\left( \begin{array}{ccc}
0 & x & 0 \\
0 & 0 & x\\
0 & 0 & 0 \end{array} \right)$ and $\tilde{Y}=\left( \begin{array}{ccc}
0 & y & 0 \\
0 & 0 & y\\
0 & 0 & 0 \end{array} \right)$.

Then $T(\tilde{X} \tilde{Y})=\left( \begin{array}{ccc}
0 & 0 & xy \\
0 & 0 & 0\\
0 & 0 & 0 \end{array} \right)\neq 0= T(\tilde{Y}) \tilde{X}^{\star}$. Hence, $T$ is not a reverse left $\star$-centralizer.
\end{ex}

Let $M_r(R)$ be a matrix ring with involution $\star$.
\begin{theorem}
Every Jordan left $\star$-centralizer over $M_r(R)$ is a reverse left $\star$-centralizer.
\end{theorem}

\begin{proof}
Let $T$ be a Jordan left $\star$-centralizer over $M_r(R)$.
\begin{align*}
&T(x^2)=T(x)x^{\star} \\
& \implies [T(x^2)]^{\star}=[T(x)x^{\star}]^{\star}=(x^{\star})^{\star}(T(x))^{\star}=x(T(x))^{\star}
\end{align*}

Now, let $S(x)=(T(x))^{\star}$ and
\begin{align*}
 S(x^2)=(T(x^2))^{\star}=x(T(x))^{\star}=xS(x).
\end{align*}

Then $S$ is a Jordan right centralizer. Hence, by Theorem \ref{lc1}, $S$ is a right centralizer. Also,
\begin{align*}
& S(xy)=xS(y) \\
& \implies (T(xy))^{\star}=x(T(y))^{\star} \\
& \implies ((T(xy))^{\star})^{\star}=(x(T(y))^{\star})^{\star} \\
& \implies T(xy)=T(y)x^{\star}.
\end{align*}

Thus, $T$ is a reverse left $\star$-centralizer over $M_r(R)$.
\end{proof}

Every left $\star$-centralizer is always a Jordan left $\star$-centralizer. The following example shows that every Jordan left $\star$-centralizer is not a left $\star$-centralizer.

\begin{ex}
We define $T:M_2(R) \rightarrow M_2(R)$ by--
\begin{center}
$T \begin{bmatrix}
x & y \\
z & t \end{bmatrix}= \begin{bmatrix}
x+y & z+t\\
0 & 0 \end{bmatrix}$.
 \end{center}

Also, we define $\star: M_2(R) \rightarrow M_2(R)$ by--\\
 \begin{center}
$\star \begin{bmatrix}
x & y \\
z & t \end{bmatrix}= \begin{bmatrix}
x & z\\
y & t \end{bmatrix}$.
 \end{center}.

 Then it is easy to see that $T$ becomes a Jordan left $\star$-centralizer.

 But $T(e_{11}e_{12})=e_{11}\neq 0=T(e_{11})e_{12}^{\star}$. Thus, $T$ is not a left $\star$-centralizer.
\end{ex}

\begin{theorem}
\label{lsc2}
If $T:M_r(R) \rightarrow M_r(R)$ is a reverse $\star$-centralizer, then there exists an $\alpha\in Z(R)$ such that $T(X)=\alpha x^{\star}$, for all $x\in M_r(R)$.
\end{theorem}

\begin{proof}
Since $T(xy)=T(y)x^{\star}$ for all $x,y \in M_r(R)$. Putting $y=1$,
\begin{equation}
T(x)=T(1)x^{\star}~~\text{for all}~~x\in M_r(R).
\end{equation}

Now, $T(xy)=y^{\star}T(x)$ for all $x,y \in M_r(R)$. Taking $y=x$ and $x=1$, we have
\begin{equation}
T(x)=x^{\star}T(1)~~\text{for all}~~x\in M_r(R).
\end{equation}

Since $\star$ is a bijective map, $T(1)\in Z(M_r(R))=Z(R)$. Therefore, $T(x)=\alpha x^{\star}$ for all $x\in M_r(R)$ where $\alpha =T(1)$.

\end{proof}

\begin{theorem}
\label{star5.3.3}
If $T:M_r(R) \rightarrow M_r(R)$ is a $\star$-centralizer, then $T=0$.
\end{theorem}

\begin{proof}
Since $T(xy)=T(x)y^{\star}=x^{\star}T(y)$ for all $x,~y \in M_r(R)$. Putting $y=1$ and $x=1$, respectively, $T(x)=x^{\star}T(1)=T(1)x^{\star}$ for all $x\in M_r(R)$.
\begin{align*}
&T(xy)=T(x)y^{\star} \implies T(1)(xy)^{\star}=T(1)x^{\star}y^{\star}
\implies T(1)[(xy)^{\star}-x^{\star}y^{\star}]=0 \\
& \implies \text{Either}~T(1)=0 ~~\text{or}~ (xy)^{\star}=x^{\star}y^{\star} \implies xy=yx.
\end{align*}

Since $xy=yx$ does not hold true in $M_r(R)$, $T(1)=0$. Thus, $T=0$.
\end{proof}

\begin{theorem}
\label{star5.3.4}
Let $m \geq 1, n \geq 1$ and $m,n \in \mathbb{Z}$, $R$ be a ring with $n(m+ n)^3$-torsion free. If $T : M_r(R) \rightarrow
M_r(R)$ be an additive mapping such that there exists a reverse $\star$-centralizer $T_0$ satisfying\\
$(m + n)T(x^2) = mx^{\star}T(x)+nT_0(x)x^{\star}$, for all $x \in M_r(R)$, then $T$ becomes reverse $\star$-centralizer. It also gives us $T=T_0$.
\end{theorem}

\begin{proof}
Let $S(x)=(T(x))^{\star}$ and $S_0(x)=(T_0(x))^{\star}$. Then $S_0(x)=\beta x,$ for some $\beta \in Z(R)$, by Theorem \ref{lsc2}. Now,
\begin{equation}
\begin{aligned}
(m+n)S(x^2)&=(m+n)(T(x^2))^{\star}=[(m+n)T(x^2)]^{\star}=[mx^{\star}T(x)+nT_0(x)x^{\star}]^{\star}\\
&=m(T(x))^{\star}x+nx(T_0(x))^{\star}=mS(x)x+nxS_0(x).
\end{aligned}
\end{equation}

Therefore, $S$ satisfies the condition of Theorem \ref{lc2}. Hence, $S$ is a two-sided centralizer, and $T$ becomes the reverse $\star$-centralizer. Since $S=S_0$, we have $T=T_0$.
\end{proof}

\begin{theorem}
Let $m \geq 1, n \geq 1$ and $m,n \in \mathbb{Z}$, $R$ be an $n(m+ n)^3$-torsion free ring. If $T : M_r(R) \rightarrow
M_r(R)$ be an additive mapping such that there exists a $\star$-centralizer $T_0$ satisfying\\
$(m + n)T(x^2) = mx^{\star}T(x)+nT_0(x)x^{\star}$, for all $x \in M_r(R)$, then $T=0$.
\end{theorem}

\begin{proof}
By Theorem \ref{star5.3.3}, $T_0=0$, $T_0$ is a reverse $\star$-centralizer. Since $T$ satisfies the conditions in Theorem \ref{star5.3.4}, we have $T=T_0=0$.
\end{proof}

Letting $S(x)=(T(x))^{\star}$ and using Theorem \ref{thm2.3}, Theorem \ref{thm2.4}, Theorem \ref{thm2.10}, Theorem \ref{thm2.11} and Theorem \ref{thm2.12}, we can prove Theorem \ref{lsc3}, Theorem \ref{th5.3.6}, Theorem \ref{th5.3.7}, Theorem \ref{th5.3.8} and Theorem \ref{lsc4}.
\begin{theorem}
\label{lsc3}
Let $R$ be a $2$-torsion free ring. If $T:M_r(R) \rightarrow M_r(R)$ is an additive map satisfying $2T(x^2)=T(x)x^{\star}+x^{\star} T(x)$, then $T$ is a reverse $\star$-centralizer and $T(x)=\alpha x^{\star}$, where $\alpha \in Z(R)$.
\end{theorem}

\begin{ex}
Let $X= \begin{bmatrix}
x & y \\
z & t \end{bmatrix}\in M_2(\mathbb{Z}_2)$ and $X^{\star}= \begin{bmatrix}
x & z \\
y & t \end{bmatrix}$.

Let us define $T:M_2(\mathbb{Z}_2)\rightarrow M_2(\mathbb{Z}_2)$ by
\begin{center}
$T(X)= \begin{bmatrix}
x+y+z+t & 0 \\
0 & x+y+z+t \end{bmatrix}$.
 \end{center}

 Then $T$ satisfies $2T(X^2)=T(X)X^{\star}+X^{\star}T(X)$, for all $X \in M_2(\mathbb{Z}_2)$.

Now, for $X=e_{11}$ and $Y= e_{12}$, $T(XY)=e_{11}+e_{12}\neq e_{11}= T(Y)X^{\star}$. Therefore, $T$ is not a reverse $\star$-centralizer. It shows that the torsion condition on $R$ is necessary for Theorem \ref{lsc3}.
\end{ex}

If $T$ satisfies the condition in Theorem \ref{lsc3}, then $T$ need not be a $\star$-centralizer.

\begin{ex}
\label{ex5.3.6}
Let us define $T(x)=x^{\star}$ for all $x\in M_2(R)$. Define $\star$ as the transpose of matrix $x$. Then $T$ satisfies $2T(x^2)=T(x)x^{\star}+x^{\star}T(x)$, for all $x \in M_2(R)$. But $T(e_{11}e_{12}) \neq T(e_{11})e_{12}^{\star}$, hence $T$ is not a $\star$-centralizer.
\end{ex}

\begin{theorem}
\label{th5.3.6}
Let $R$ be a $2$-torsion free ring. If $T:M_r(R) \rightarrow M_r(R)$ is an additive map satisfying $2T(x^2)=T_0(x)x^{\star}+x^{\star} T(x)$ where $T_0:M_r(R)\rightarrow M_r(R)$ is an additive map satisfying $2T_0(x^2)=T_0(x)x^{\star}+x^{\star} T_0(x)$, then $T$ is a reverse $\star$-centralizer and $T(x)=\alpha x^{\star}$, where $\alpha \in Z(R)$.
\end{theorem}

\begin{theorem}
\label{th5.3.7}
Let $R$ be a ring. If $T$ is an additive map on $M_r(R)$ with $T(xyx) = x^{\star}T(y)x^{\star}$, for all $x, y \in M_r(R)$, then $T$ becomes reverse $\star$-centralizer. In particular, $T(x)=\alpha x^{\star}$, where $\alpha \in Z(R)$.
\end{theorem}

\begin{theorem}
\label{th5.3.8}
Let $R$ be a ring. If $T$ is an additive map on $M_r(R)$ with $2T(xyx) = T(x)y^{\star}x^{\star} + x^{\star}y^{\star}T(x)$ for all $x, y \in M_r(R),$ then $T$ is a reverse $\star$-centralizer. In particular, $T(x)=\alpha x^{\star}$ where $\alpha \in Z(R)$.
\end{theorem}

\begin{theorem}
\label{lsc4}
Let $R$ be a ring with $2$-torsion free. If $T$ is an additive map on $M_r(R)$ with
$3T(xyx) = T(x)y^{\star}x^{\star} + x^{\star}T(y)x^{\star} + x^{\star}y^{\star}T(x)$ for all $x, y \in M_r(R)$,
then $T$ is a reverse $\star$-centralizer. In particular, $T(x)=\alpha x^{\star}$, where $\alpha \in Z(R)$.
\end{theorem}

\begin{ex}
Let $A= \begin{bmatrix}
a & b \\
c & d \end{bmatrix}\in M_2(\mathbb{Z}_2)$ and $A^{\star}= \begin{bmatrix}
a & c \\
b & d \end{bmatrix}$.

Let $T:M_2(\mathbb{Z}_2)\rightarrow M_2(\mathbb{Z}_2)$ be defined by
\begin{center}
$T(A)= \begin{bmatrix}
0 & c \\
b & 0 \end{bmatrix}$.
 \end{center}

 Then $T$ satisfies $3T(ABC) = T(A)B^{\star}A^{\star} + A^{\star}T(B)A^{\star} + A^{\star}B^{\star}T(A)$, for all $A,~B \in M_2(\mathbb{Z}_2)$.

Now, for $A=e_{11}$ and $B= e_{12}$, $T(AB)=e_{21} \neq 0=B^{\star} T(A)$. Therefore, $T$ is not a reverse $\star$-centralizer. It shows that the torsion condition on $R$ is necessary for the Theorem \ref{lsc4}.
\end{ex}

Example \ref{ex5.3.6} shows that even $T$ satisfies the conditions in
Theorems \ref{th5.3.6}, \ref{th5.3.7}, \ref{th5.3.8} and \ref{lsc4}, but need not be a $\star$-centralizer.
Theorem \ref{lc2} motivates us to post a conjecture as follows.
\begin{conj}
Let $m \geq 1$, $n \geq 1$ and  $m,~n\in \mathbb{Z}$, $R$ be a ring with some suitable torsion restrictions. If $T : M_r(R) \rightarrow M_r(R)$ be an additive mapping such that
\begin{equation}
(m  + n)T(x^2) = mT(x)x + nxT(x),~\text{for all}~ x \in M_r(R),
\end{equation}
then $T$ become a two-sided centralizer.
\end{conj}

\section{ Conclusions}
In this paper, we discussed centralizers over matrix rings. Then we proved that the Jordan left centralizer over the matrix ring becomes the left centralizer, which is generally not true. Also, we showed that every two-sided centralizer over the matrix ring is of some particular form. Later, it showed that the map becomes a two-sided centralizer for the matrix ring when it satisfies some functional equation. We also prove every Jordan left $\star$-centralizer is a reverse left $\star$-centralizer over matrix ring with involution. Finally, we have established that every reverse $\star$-centralizer over $\star$-matrix ring has some particular form, and a map satisfying some equations over $\star$-matrix ring becomes reverse $\star$-centralizer.

\section{Acknowledgement}

The authors thank the Department of Science and Technology, Govt. of India, for financial support under DST/INSPIRE Fellowship/IF140850 and the Indian Institute of Technology Patna for providing the research facilities. The authors would also like to thank the anonymous referee(s) and the Editor-in-Chief for their valuable suggestions to improve the presentation of the manuscript.


\begin{thebibliography}{99}

\bibitem{ali2012generalized} Ali, S., (2012), On generalized $\star$-derivations in $\star$-rings, Palest. J. Math., 1(1), pp. 32-37.

\bibitem{shakir2010jordan} Ali, S., Fo\v{s}ner, A., (2010), On {Jordan} ($\alpha$, $\beta$)$^{\star}$-Derivations in Semiprime $\star$-Rings, Int. J. Algebra, 4(3), pp. 99-108.


\bibitem {ali2014generalized} Ali, S., Fo\v{s}ner, A., (2014), On generalized $(m, n)$-derivations and generalized $(m, n)$-Jordan derivations in rings, Algebr. Colloq., 21(3), pp. 411-420.



\bibitem {ali2013jordan} Ali, S., Dar, N. A., Vukman, J., (2013), Jordan left $\star$-centralizers of prime and semiprime rings with involution, Beitr{\"a}ge Algebra Geom., 54(2), pp. 609-624.



\bibitem {beidar1995rings} Beidar, K. I.,  Martindale III, W. S., Mikhalev, A. V., (1996), Rings with Generalized Identities, Marcel Dekker, Inc. New York.

\bibitem{brevsar1992structure} Bre{\v{s}}ar, M., Zalar, B., (1992), On the structure of Jordan $\star$-derivations, Colloq. Math., 63(2), pp. 163-171.

\bibitem{dar2020structure} Dar, N. A., Ali, S., (2021), On the structure of generalized Jordan $\star$-derivations of prime rings, Comm. Algebra, 49(4), pp. 1422-1430.

\bibitem {fovsner2013note} Fo\v{s}ner, A., (2013), A note on generalized $(m, n)$-Jordan centralizers, Demonstratio Math., 46(2), pp. 254-262.

\bibitem{ghosh2019new} Ghosh, A., Prakash, O., (2019), New Results on Generalized $(m, n)$-Jordan Derivations over Semiprime Rings, Southeast Asian Bull. Math., 43(3), pp. 323-331.

\bibitem{johnson1964introduction} Johnson, B. E., (1964), An introduction to the theory of centralizers, Proc. London Math. Soc., 3(2), pp. 299-320.



\bibitem{vukman1999identity} Vukman, J., (1999), An identity related to centralizers in semiprime rings, Comment. Math. Univ. Carolin., 40(3), pp. 447-456.

\bibitem {vukman2001centralizers} Vukman, J., (2001), Centralizers of semiprime rings, Comment. Math. Univ. Carol., 42(2), pp. 237-245.

\bibitem{vukman2003centralizers} Vukman, J.,  Kosi-Ulbl, I., (2003), On centralizers of semiprime rings, Aequat. Math., 66(3), pp. 277-283.

\bibitem{vukman2003equation} Vukman, J.,  Kosi-Ulbl, I., (2003), An equation related to centralizers in semiprime rings, Glas. Mat., 38(2), pp. 253-261.

\bibitem{vukman2008m} Vukman, J., (2008), On $(m, n)$-Jordan derivations and commutativity of prime rings, Demonstratio Math., 41(4), pp. 773-778.

\bibitem {vukman2010m} Vukman, J., (2010), On $(m, n)$-Jordan centralizers in rings and algebras, Glas. Mat., 45(1),  pp. 43-53.

\bibitem{wendel1952left}  Wendel, J. G., (1952), Left centralizers and isomorphisms of group algebras, Pac. J. Math., 2(2), pp. 251-261.

\bibitem{zalar1991centralizers}  Zalar, B., (1991) On centralizers of semiprime rings, Comment. Math. Univ. Carolin, 32(4), pp. 609-614.

\end{thebibliography}
\end{document}